\documentclass[10pt]{amsart}
\usepackage{graphicx,amssymb,amsfonts,amsmath,amsthm,newlfont,dsfont}
\usepackage[foot]{amsaddr}
\usepackage{epsfig,url}
\usepackage{color}
\usepackage{bm}

\newlength\bshft
\def\fakebold#1{\setbox0=\hbox{$#1$}#1\kern-\wd0\kern\bshft#1\kern-\wd0\kern\bshft#1}

\usepackage[sort&compress,numbers]{natbib}

\usepackage[all,2cell]{xy} \UseAllTwocells \SilentMatrices

\newtheorem{thm}{Theorem}[section]

\newtheorem{prop}[thm]{Proposition}
\newtheorem{conj}[thm]{Conjecture} 
\theoremstyle{definition}
\newtheorem{defi}[thm]{Definition}

\newtheorem{example}{Example}

\theoremstyle{remark}

\numberwithin{equation}{section}
\newcommand{\Z}{\mathbb Z}
\newcommand{\C}{\mathbb C}

\newcommand{\R}{\mathbb R}

\newcommand{\Q}{\mathbb Q}

\newcommand{\SL}{\mathrm{SL}}

\newcommand{\Res}{\mathrm{Res}}

\newcommand{\HH}{\mathfrak{H} }

\newcommand{\cM}{\mathcal{M} }

\newcommand{\3}{\mathsf{3} }
\newcommand{\5}{\mathsf{5} }

\newcommand{\matix}[1]{\left(\begin{smallmatrix}#1 \end{smallmatrix}\right)}

\def\beq{\begin{equation}}
\def\eeq{\end{equation}}
\def\bsp#1\esp{\begin{split}#1\end{split}}

\newcommand{\mati}[4]{\left(\begin{smallmatrix} #1 & #2 \\ #3 & #4\end{smallmatrix}\right)}

\newcommand{\cMmag}{\cM^{\mathrm{mag}}}
\newcommand{\cMsm}{\cM^{\mathrm{s.m.}}}
\newcommand{\rd}{\textrm{d}}
\newcommand{\cF}{\mathcal{F}}

\renewcommand{\Im}{\textrm{Im}}
\begin{document}
\begin{flushright}BONN-TH-2024-06\end{flushright}

\author{Kilian B\"onisch$^{\textrm{A,B}}$}
\address[A]{Max-Planck-Institut f\"ur Mathematik, Bonn, D-53111, Germany}
\author{Claude Duhr$^{\textrm{B}}$}

\author{Sara Maggio$^{\textrm{B}}$}
\address[B]{Bethe Center for Theoretical Physics, Universit\"at Bonn, D-53115, Germany}

\email{boenisch@physik.uni-bonn.de}
\email{cduhr@uni-bonn.de}
\email{smaggio@uni-bonn.de}

\begin{title}[Some conjectures around magnetic modular forms]
{Some conjectures around magnetic modular forms} \end{title}
\maketitle

\begin{abstract}
We study a class of meromorphic modular forms characterised by Fourier coefficients that satisfy certain divisibility properties. We present new candidates for these so-called magnetic modular forms, and we conjecture properties that these functions should obey. In particular, we conjecture that magnetic modular forms are closed under the standard operators acting on spaces of modular forms ($\SL_2(\Z)$ action, Hecke and Atkin-Lehner operators), and that they are characterised by algebraic residues and vanishing period polynomials. We use our conjectures to construct examples of real-analytic modular forms with poles. 
%\kiliancomment{I am still quite unhappy with ``real-analytic with poles'', but I don't have a better snappy description.}
\end{abstract}

% !TEX root = main.tex

\section{Introduction}
\label{sec:intro}

Holomorphic modular forms for finite-index subgroups of the modular group $\SL_2(\Z)$ have been studied by mathematicians for over a century. They have many interesting arithmetic properties, and are the foundation for many deep results in modern mathematics. More recently, it has become clear that modular forms also play an important role in the computation of Feynman integrals that arise in perturbative Quantum Field Theory, cf. ref.~\cite{Adams:2017ejb}.

Meromorphic modular forms, which are allowed to have poles in the upper half plane and at the cusps, are much less studied, and their structure and their arithmetic properties are much less understood. There is an interesting subclass of meromorphic modular forms, whose coefficients of the Fourier expansion satisfy certain congruence relations (see Section~\ref{sec:definitions} for the precise definition). This subclass was called \emph{magnetic} in ref.~\cite{magnetic2}, because the first non-trivial candidate was constructed in ref.~\cite{magnetic1} inspired by a physics problem (and magneticity was later proven in ref.~\cite{magnetic2}). Additional examples were subsequently given in ref.~\cite{magnetic3} using the Shimura-Borcherds lift. However, a general theory and an understanding of the algebraic and arithmetic properties of magnetic modular forms is still lacking.

Recently, new candidates for magnetic modular forms for certain congruence subgroups have been constructed from two independent sources. First, in his PhD thesis \cite{ThesisBoenisch}, which is based on the work in progress \cite{FiberingOutCalabiYauMotives}, one the authors has studied meromorphic modular forms associated with degenerations of families of algebraic varieties. These seem to be part of extensions of motives of holomorphic Hecke eigenforms and computations suggest that they are Hecke eigenforms modulo magnetic modular forms. Second, magnetic modular forms seem to arise in Quantum Field Theory when considering differential equations satisfied by dimensionally-regulated Feynman integrals related to certain one-parameter families of K3 surfaces~\cite{Pogel:2022yat,epsilon_form}. While in both cases it is not proven that these modular forms are indeed magnetic, it was checked that several hundred first  coefficients of the Fourier expansion have the desired divisibility properties. 
 This goal of this paper is to use the available candidates for magnetic modular forms to distill for the first time a list of algebraic and arithmetic properties that this class of modular forms should obey, and to propose a set of conjectures that characterise magnetic modular forms.
 %Remarkably, we find that all known candidates are closed under the standard operators acting on spaces of modular forms ($\SL_2(\Z)$ action, Heck and Atkin-Lehner operators). Moreover, all candidates have very simple residues, and vanishing period polynomials.
% !TEX root = main.tex

\section{Definitions}
\label{sec:definitions}

Let $k$ be some integer and let $\Gamma$ be a finite-index subgroup of $\SL_2(\Z)$. The main object of interest in this paper is the complex vector space $\cM_k(\Gamma)$ of meromorphic modular forms of {weight} $k$ under $\Gamma$. This consists of all meromorphic functions $f$ on the upper half-plane $\HH = \{\tau\in\C \ | \ \Im \, \tau>0\}$ that satisfy $f|_{k}\gamma=f$ for all $\gamma\in \Gamma$ and that are meromorphic at the cusps. Here we use the usual weight $k$ action of $\gamma = \mati{a}{b}{c}{d} \in \text{GL}_2^+(\R)$ given by
\beq
(f|_{k}\gamma)(\tau) \, = \, (c \, \tau+d)^{-k}\,f\left(\frac{a \, \tau+b}{c \, \tau+d}\right)\, .
\eeq
Note that $f$ is allowed to have poles in $\HH$ and at the cusps. Important finite-dimensional subspaces of $\cM_k(\Gamma)$ are the space of holomorphic modular forms $M_k(\Gamma)$ and the space of cusp forms $S_{k}(\Gamma)$. 
In the following we will only consider the case where $\Gamma$ is a \emph{congruence subgroup} of some \emph{level} $N \in \mathbb{N}$, i.e., it contains the subgroup
\beq
\Gamma(N) \, = \, \{\gamma\in \SL_2(\Z) \ | \ \gamma \equiv \mathds{1}\!\!\!\!\mod N\}\, .
\eeq
Important examples of congruence subgroups are
\begin{align}
    \Gamma_1(N) \, &= \, \{\mati{a}{b}{c}{d}\in \SL_2(\Z) \ | \ a,d \equiv 1\!\!\!\!\mod N, \ c \equiv 0\!\!\!\!\mod N \}\,, \\
    \Gamma_0(N) \, &= \, \{\mati{a}{b}{c}{d}\in \SL_2(\Z) \ | \ c \equiv 0\!\!\!\!\mod N\}\, .
\end{align}

\begin{defi}
\label{def:magnetic}
Let $f$ be a meromorphic modular form of level $N$ with Fourier expansion
\beq
f(\tau) \, = \, \sum_{n \in \mathbb{Z} \setminus \{0\}} a_n\,q^{n/N}\,,\qquad a_n \, \in \, \overline{\Q}\,,\qquad q \, = \, e^{2\pi i\tau}\,.
\eeq
We say that $f$ is \emph{magnetic of depth} $d > 0$ if the algebraic numbers $a_n/n^d$ have globally bounded denominators.
\end{defi}

In essence, a meromorphic modular form is magnetic of depth $d$ if the $n^{\textrm{th}}$ coefficient in the Fourier expansion is divisible by $n^d$. A folklore conjecture states that there are no (non-zero) holomorphic magnetic modular forms.

Our definition of magnetic modular forms is slightly more general than the one originally presented in ref.~\cite{magnetic2}, because we consider divisibility in $\overline{\Q}$ and because we allow bounded denominators. This permits us to define the $\overline{\Q}$-vector space $\cMmag_{k,d}(\Gamma)$ of modular forms of weight $k$ and depth $d$. There is an obvious filtration by depth 
\beq
\cM_k(\Gamma) \, \supseteq \, \cMmag_{k,1}(\Gamma) \, \supseteq \, \cMmag_{k,2}(\Gamma) \, \supseteq \, \ldots
\eeq
and we also write $\cMmag_{k}(\Gamma) = \bigcup_{d\ge 1}\cMmag_{k,d}(\Gamma)$.

One example of a magnetic modular form $\phi \in \cMmag_{4,1}(\Gamma_0(8))$ is given by 
\beq
\phi \, = \, - \frac{1-6 \, t^2+t^4}{(1+t^2)^2} \, f\,,
\eeq
where $t$ is the Hauptmodul of $\Gamma_0(8)$ that is normalised by $t(\tau) = 1-8 \, q+\ldots$ and $f$ is the unique newform in $S_4(\Gamma_0(8))$. These can be represented by
\beq\label{eq:Gamma0(8)}
t(\tau) \, = \, \frac{\eta(\tau)^8 \, \eta(4\, \tau)^4}{\eta(2 \, \tau)^{12}} \qquad \text{and} \qquad f(\tau) \, = \, \eta(2 \, \tau)^4 \, \eta(4 \, \tau)^4\,,
\eeq
with the Dedekind eta function
\beq
\eta(\tau) = e^{i\pi\tau/12}\prod_{n=1}^\infty \big(1-e^{2\pi in\tau}\big)\,.
\eeq
The divisibility of the first Fourier coefficients can be observed in the expansion
\begin{align}
    \phi(\tau) \, = \, q - 132 \, q^3 + 5630 \, q^5 - 189672 \, q^7 + 5768181 \, q^9 + \cdots \, .
\end{align}
This example originates from a computation in physics and coined the term ``magnetic''. $\phi$ was conjectured to be magnetic in ref.~\cite{magnetic1} and this was proven in ref.~\cite{magnetic2} by expressing $\phi$ as the Borcherds-Shimura lift of a modular form of half-integral weight.

Other examples of magnetic modular forms for the full modular group $\SL_2(\Z)$ are given in ref.~\cite{magnetic3}. Two examples with $(k,d) = (4,1)$ and $(k,d) = (6,2)$ are given by
\begin{align}
    \frac{E_4(\tau)}{j(\tau)} \, &= \, q - 504 \, q^2 + 180252 \, q^3 - 56364992 \, q^4 + 16415391870 \, q^5 + \cdots \\
    \frac{E_6(\tau)}{j(\tau)} \, &= \, q - 1248 \, q^2 + 714996 \, q^3 - 307862528 \, q^4 + 114237828150 \, q^5 + \cdots 
\end{align}
in terms of Eisenstein series and the $j$-invariant.

\section{Operators on spaces of magnetic modular forms}
\label{sec:sl2_action}

There are important operators on the spaces $\mathcal{M}_k(\Gamma(N))$ and a natural question is whether these preserve magneticity. 

\subsection{Action by $\boldsymbol{\SL_2(}\protect\fakebold{\mathbb{Z}}\boldsymbol{)}$}

The group $\Gamma(N)$ is normal in $\SL_2(\Z)$, and hence we have a natural action of $\Gamma(N) \backslash \SL_2(\Z)$ on $\cM_k(\Gamma(N))$. We have computed this action on many examples of magnetic modular forms (see Appendix \ref{sec:NumericalSL2Action}) and our computations suggest the following.

\begin{conj}\label{conj:action}
For any $\gamma \in \SL_2(\Z)$, the action $f \mapsto f|_k\gamma$ on $\cM_k(\Gamma(N))$ can be restricted to the magnetic subspaces $\cM^{\mathrm{mag}}_{k,d}(\Gamma(N))$.
\end{conj}

Since all cusps are $\SL_2(\Z)$-equivalent, the conjecture says that magneticity does not depend on the chosen cusp for the Fourier expansion. 

For example, acting with the twelve cosets of $\Gamma_0(8)$ in $\SL_2(\Z)$ on the magnetic modular form $\phi$ one obtains the magnetic modular forms
\begin{align*}
    \pm\phi(\tau) \, , \ \pm \tfrac{1}{4} \, \phi(\tfrac{\tau}{2}+\tfrac{1}{4}) \, , \ \pm \tfrac{1}{64} \, \phi(\tfrac{\tau}{8}) \, , \ \pm \tfrac{1}{64} \, \phi(\tfrac{\tau}{8}+\tfrac{1}{8}) \, , \ \pm \tfrac{1}{64} \, \phi(\tfrac{\tau}{8}+\tfrac{2}{8}) \, , \ \pm \tfrac{1}{64} \, \phi(\tfrac{\tau}{8}+\tfrac{3}{8})\,.
\end{align*}

\subsection{Action by Hecke operators}

The map $f(\cdot) \mapsto f(N \, \cdot)$ gives an isomorphism 
\beq\bsp
    \mathcal{M}_k(\Gamma(N)) \, &\cong \, \bigoplus_{\chi \text{ mod } N} \mathcal{M}_k(\Gamma_0(N^2),\chi) \, ,
\esp\eeq
where $\mathcal{M}_k(\Gamma_0(N^2),\chi)$ denotes the subspace of $\mathcal{M}_k(\Gamma_1(N^2))$ that transforms under $\Gamma_1(N^2) \backslash \Gamma_0(N^2)$ according to a character $\chi$. We can identify the character $\chi$ with a Dirichlet character by using the isomorphism $\Gamma_1(N^2) \backslash \Gamma_0(N^2) \cong (\mathbb{Z}/N^2 \mathbb{Z})^\times$ given by ${\mati{a}{b}{c}{d} \mapsto a}$, and the sum above is only over Dirichlet characters of modulus $N$. The spaces $\mathcal{M}_k(\Gamma_0(N),\chi)$ are equipped with an action by Hecke operators. For any $n \in \mathbb{N}$ coprime to $N$, the action of the Hecke operator $T_n$ on $\cM_k(\Gamma_0(N),\chi)$ is defined by 
\begin{align}
    T_n(f)\, = \, n^{k-1} \, \sum_{\substack{a \cdot d = n\\ 0 \leq b < d}} \chi(d) \, f|_k \mati{a}{b}{0}{d} \, . 
\end{align}
This gives a simple action on the Fourier coefficients and it is easy to see that the Hecke operators preserve $\cM_k(\Gamma_0(N),\chi)$ for depths $d \leq k-1$. This proves the following proposition.

\begin{prop}
The action of Hecke operators on $\cM_k(\Gamma_0(N),\chi)$ can be restricted to the magnetic subspaces $\cM^{\mathrm{mag}}_{k,d}(\Gamma_0(N),\chi)$ for depths $d \leq k-1$.
\end{prop}

It follows that one can use Hecke operators to construct infinitely many linearly independent magnetic modular forms. For example, acting on $E_4/j$ with $T_2$ and $T_3$ one obtains the magnetic modular forms
\begin{align}
    T_2\left(\frac{E_4}{j}\right) \, &= \, -72 \, \frac{30375+7 \, j}{(54000-j)^2} \, E_4\,, \\ 
    T_3\left(\frac{E_4}{j}\right) \, &= \, 108 \, \frac{4194304000000 - 1339392000 \, j+1669 \, j^2}{j \,  (12288000+j)^2} \, E_4  \, .
\end{align}

\subsection{Action by Atkin-Lehner involutions}

The map $f(\cdot) \mapsto f(N \, \cdot)$ allows us to identify $\mathcal{M}_k(\Gamma(N))$ with a subspace of $\mathcal{M}_k(\Gamma_1(N^2))$. The spaces $\mathcal{M}_k(\Gamma_1(N))$ are equipped with an action by Atkin-Lehner involutions. For any exact divisor $Q$ of $N$ (i.e.\ $Q | N$ and $(Q,N/Q)=1$), the Atkin-Lehner involution $W_Q$ on $\cM_k(\Gamma_1(N))$ is defined by
\begin{align*}
    W_Q(f) \, = \,  f|_k \tfrac{1}{\sqrt{Q}} \mati{Q \, x}{y}{N \, z}{Q \, w}  \, ,
\end{align*}
where $x,y,z,w$ are any choice of integers satisfying $\det \mati{Q \, x}{y}{N \, z}{Q \, w} = Q$ and the congruences $y \equiv 1 \text{ mod } Q$ and $x \equiv 1 \text{ mod } N/Q$. We have computed the action of Atkin-Lehner involutions on several examples of magnetic modular forms (see Appendix \ref{sec:NumericalALAction}) and our computations suggest the following:

\begin{conj}
The action of Atkin-Lehner involutions on $\cM_k(\Gamma_1(N))$ can be restricted to the magnetic subspaces $\cM^{\mathrm{mag}}_{k,d}(\Gamma_1(N))$.
\end{conj}

For example, the magnetic modular form $\phi$ satisfies $W_8(\phi) = - \phi$.

\section{Poles and residues of magnetic modular forms}
\label{sec:poles}

\subsection{Strongly magnetic modular forms}
\label{sec:strongly_magnetic}

One can construct trivial examples of magnetic modular forms by using the derivative $D = \frac{1}{2\pi i}\partial_{\tau} = q \, \partial_q$. For any weight $k \geq 2$, Bol's identity~\cite{Bol} implies that the $(k-1)^{\textrm{th}}$ derivative of a modular form of weight $2-k$ is a modular form of weight $k$, i.e.\ we have a map
\beq
D^{k-1} : \cM_{2-k}(\Gamma) \to \cM_{k}(\Gamma)\,.
\eeq
If we restrict to the subspace $\cM_{2-k}^{\text{alg}}(\Gamma)$ with algebraic Fourier coefficients (which necessarily have bounded denominators since $\Gamma$ is a congruence subgroups), we obtain a map 
\beq
D^{k-1} : \cM_{2-k}^{\text{alg}}(\Gamma) \to \cM_{k,k-1}^{\text{mag}}(\Gamma)\,.
\eeq
Clearly, we are interested in those magnetic modular forms that do not lie in the image of $D^{k-1}$. Note that, away from the cusps, the poles of elements of $D^{k-1}\cM_{2-k}^{\textrm{alg}}(\Gamma)$ have order at least $k$.

\begin{defi}
A strongly magnetic modular form of even weight $k>2$ is a magnetic modular form of weight $k$ and depth $k/2-1$ whose poles all have order less or equal than $k/2$. 
\end{defi}

We denote the vector space of strongly magnetic modular forms by $\cMsm_k(\Gamma)$.
All examples of magnetic modular forms we are aware of are linear combinations of strongly magnetic modular forms and elements in the image of $D^{k-1}$ (see Appendix~\ref{app:candidates}). Based on this observation, we propose the following:

%This is in particular the case for all examples in refs.~\cite{magnetic1,magnetic2,magnetic3}, and also for all candidates presented in section~\ref{sec:candidates}. In particular, there are no known candidates for magnetic modular forms of odd weights. In particular, we have used the procedure of appendix~\ref{app:congruence} to search for candidates of weight 3 and 5. While this procedure allows us to recover the candidates of even weight presented in ~\ref{sec:candidates}, it fails to find any candidate in odd weights. 

%The function $C_{6a}$ is a magnetic modular form of weight 6 and depth 2, and it has a triple pole at $\tau\in\Gamma_0(2)\tau_0$, with $\tau_0=\frac{i}{\sqrt{2}}$ and $A(\tau_0)=\frac{1}{64}$. The functions $C_4$ and $C_{6b}$, however, seem to have a pole at $\tau\in\Gamma_0(2)\tau_1$, with $\tau_1=\frac{i-1}{2}$ and $A(\tau_1)=-\frac{1}{64}$, of order 1 and 2 respectively. Hence, at first glance it appears that $C_4$ and $C_{6b}$ are magnetic, but not strongly magnetic. However, $\tau_1$ is an elliptic point for $\Gamma_0(2)$, and if we pass to an appropriate cover where $\tau_1$ is not elliptic (for example $\Gamma(2)$), we see that $C_4$ and $C_{6b}$ have a pole of order 2 and 3 respectively, and so they are strongly magnetic. Based on these observations, we propose:

\begin{conj}
\label{conj:pole}
There is a decomposition
\beq
\cMmag_{k}(\Gamma) = \left\{\begin{array}{ll}
\cMsm_{k}(\Gamma) \oplus D^{k-1}\cM_{2-k}^{\mathrm{alg}}(\Gamma)\,,& \textrm{~~$k$ even}\,,\\
D^{k-1}\cM_{2-k}^{\mathrm{alg}}(\Gamma)\,,& \textrm{~~$k$ odd}\,.
\end{array}\right.
\eeq
\end{conj}

\subsection{Residues of magnetic modular forms}
% Let $f$ be a non-zero strongly magnetic modular form of weight $k$. Conjecture~\ref{conj:pole} implies that it has even weight and it has only poles of order exactly $\frac{k}{2}$. 

We have computed the residues of many examples and candidates of magnetic modular forms (see Appendix \ref{sec:NumericalResidues}). In all cases, we find that the residues take a particularly simple form, which we expect to hold in general.

%In particular, we have computed the residues of all candidates in section~\ref{sec:candidates} and of all examples presented in refs.~\cite{magnetic1,magnetic2,magnetic3}. 

\begin{conj}
\label{conj:residues}
Let $f$ be a magnetic modular form of weight $k$. Then the residues of $\tau^m \, f(\tau)$ for $0\le m\le k-2$ are in $\frac{1}{(2\pi i)^{k/2}} \overline{\mathbb{Q}}$ and only non-vanishing at complex multiplication points.
\end{conj}

\begin{example}\label{ex:phi_residues}
Consider the magnetic modular form $\phi$. It has poles only in the $\Gamma_0(8)$-orbits of $\tau_{\pm} = \pm \tfrac{1}{4}+\tfrac{i}{4}$, and the residues are given by 
\begin{align}
    \text{Res}_{\tau_\pm} \ \begin{pmatrix} 1 \\ \tau \\ \tau^2\end{pmatrix} \phi(\tau)\ \, = \, \frac{1}{(2\pi i)^2} \begin{pmatrix} \mp 1 \\ -1/4 \\ \mp 1/8\end{pmatrix} \, .
\end{align}
\end{example}

\begin{example}\label{ex:SL2Ej}
The magnetic modular forms $E_4/j$ and $E_6/j$ have poles only in the $\SL_2(\Z)$-orbit of $\tau_0 = e^{2\pi i / 3} = -\tfrac{1}{2}+\tfrac{\sqrt{-3}}{2}$ and the residues are
\begin{align*}
    \text{Res}_{\tau_0} \left( \begin{array}{c}
         1\\
         \tau\\
         \tau^2
    \end{array}\right) \frac{E_4(\tau)}{j(\tau)} \, &= \, \frac{\sqrt{-3}}{(2\pi i)^2} \left( \begin{array}{c}
         -1/288  \\
         1/576 \\
         -1/288
    \end{array}\right)\,, \\
    \text{Res}_{\tau_0} \left( \begin{array}{c}
         1\\
         \tau\\
         \tau^2\\
         \tau^3\\
         \tau^4
    \end{array}\right) \frac{E_6(\tau)}{j(\tau)} \, &= \, \frac{1}{(2\pi i)^3} \left( \begin{array}{c}
         -1/32  \\
         1/64 \\
         -1/64 \\
         1/64 \\
         -1/32
    \end{array}\right) \, .
\end{align*}
\end{example}

% !TEX root = main.tex

\section{Periods of magnetic modular forms}
\label{sec:periods}

For any ring $R$ we denote by $V_{k}(R)$ the $R$-module of homogeneous polynomials of degree $k$ in two variables $X$ and $Y$. Elements $\gamma = \mati{a}{b}{c}{d}\in\SL_2(\Z)$ act on $V_k(R)$ by
\beq\label{eq:pol_action}
P(X,Y)|\gamma = P(a \, X+b \, Y,c \, X+d \, Y)\, .
\eeq
%If $f$ is a meromorphic modular form of weight $k$ for $\Gamma\subseteq \SL_2(\Z)$, we can define the differential one-form
%\beq
%\underline{f} := \rd \tau\,f(\tau)\,(X-\tau Y)^{k-2}\,.
%\eeq
%It is easy to check that this differential form is $\Gamma$-invariant:
%\beq
%\Big(\underline{f}_{|k\gamma}\Big)|_{\gamma} = \underline{f}\,.
%\eeq

For any weight $k \geq 2$, let $\mathcal{S}_k(\Gamma)$ be the space of meromorphic modular forms $f \in \mathcal{M}_k(\Gamma)$ whose residues $\text{Res}_{\tau_0} \, \tau^m \, f(\tau)$ for $0\le m\le k-2$ vanish at all $\tau_0 \in \HH$. For any $f \in \mathcal{S}_k(\Gamma)$ and any $\tau_0 \in \overline{\HH}$ at which $f$ is holomorphic\footnote{Here, $f$ is said to be holomorphic at a cusp if it vanishes at the cusp.} we can define the map  
\begin{align}\label{eq:cocyle_def}
C_{f,\tau_0} : \ \Gamma \, &\rightarrow \, V_{k-2}(\mathbb{C})\\
\nonumber\gamma \, &\mapsto \, \frac{(2\pi i)^{k-1}}{(k-2)!} \int_{\gamma^{-1}\tau_0}^{\tau_0} (X-\tau \, Y)^{k-2}\, f(\tau) \, \rd\tau\,.
\end{align}
It is easy to check that $C_{f,\tau_0}$ satisfies the cocycle condition
\beq\label{eq:cocycle}
C_{f,\tau_0}(\gamma\gamma') = C_{f,\tau_0}(\gamma)|\gamma' + C_{f,\tau_0}(\gamma')\,.
\eeq
While $C_{f,\tau_0}$ depends on the choice of $\tau_0$, the associated class $C_f = [C_{f,\tau_0}]$ in the first group cohomology 
\beq
H^1(\Gamma,V_{k-2}(\mathbb{C})) = \frac{\big\{C : \Gamma\to V_{k-2}(\mathbb{C}) \ | \ d^{(2)} C \, = \, 0 \big\}}{\big\{d^{(1)}P  \, | \, P\in V_{k-2}(\mathbb{C})\big\}}
\eeq
does not. Here, we use the coboundary operators defined by
\begin{align*}
    (d^{(1)} P)(\gamma) \, &= \, P|\gamma - P\,, \\
    (d^{(2)} C)(\gamma,\gamma') \, &= \, C(\gamma)|\gamma' + C(\gamma') - C(\gamma \gamma') \, .
\end{align*}
In fact, it was shown by Eichler \cite{EichlerAbelscheIntegrale} that one obtains an isomorphism
\begin{align}
    \mathcal{S}_k(\Gamma) / D^{k-1} \mathcal{M}_{2-k}(\Gamma) \, &\rightarrow \, H^1(\Gamma,V_{k-2}(\mathbb{C})) \\
 \nonumber   [f] \, &\mapsto \, C_f \, .
\end{align}

We now generalise the previous construction to elements $f \in \mathcal{M}_k(\Gamma)$ which have non-vanishing residues. Let $\Lambda_f \subset \mathbb{C}$ be the $\mathbb{Q}$-vector space spanned by the residues $(2\pi i)^k \, \text{Res}_{\tau_0} \, \tau^m \, f(\tau)$ for $0\le m\le k-2$ and $\tau_0 \in \HH$. For any choice of $\tau_0 \in \overline{\HH}$ at which $f$ is holomorphic, we can again consider the integral 
\begin{align}
    \frac{(2\pi i)^{k-1}}{(k-2)!} \int_{\gamma^{-1}\tau_0}^{\tau_0} (X-\tau \, Y)^{k-2}\, f(\tau) \, \rd\tau \, ,
\end{align}
but this now depends on the path of integration. However, taking the quotient by contributions from residues, we obtain a well defined map $C_{f,\tau_0} : \Gamma \rightarrow V_{k-2}(\mathbb{C}/\Lambda_f)$ and an associated cohomology class $C_f \in H^1(\Gamma,V_{k-2}(\mathbb{C}/\Lambda_f))$. Based on numerical observations (see Appendix \ref{sec:NumericalPeriods}) we conjecture the following:

\begin{conj}\label{conj:periods}
Let $f \in \mathcal{M}_k(\Gamma)$ be a meromorphic modular form of some weight $k>2$ and with algebraic Fourier coefficients. Further suppose that the residues of $f$ are algebraic, by which we mean that $\Lambda_f \subset (2\pi i)^{k/2} \overline{\mathbb{Q}}$. Then $f$ is magnetic if and only if $C_f = 0$.
\end{conj}

Note that the statement of the conjecture is not true for weight $k=2$. For example, for $f = \frac{1}{2\pi i} \frac{j'}{j} = \frac{1}{2\pi i} \log(j)'$ it is clear that for any $\gamma \in \SL_2(\Z)$ and any path of integration we have 
\begin{align}
2 \pi i \int_{\gamma^{-1}\tau_0}^{\tau_0} f(\tau) \, \rd \tau \, \in \, 2\pi i \, \mathbb{Z}
\end{align}
but $f$ is not magnetic. 

Also note that for $\Lambda_f \subset (2\pi i)^{k/2} \overline{\mathbb{Q}}$ we have an isomorphism 
\begin{align}
    H^1(\Gamma,V_{k-2}(\mathbb{C}/\Lambda_f)) \, \cong \, H^1(\Gamma,V_{k-2}(\mathbb{Q})) \otimes_{\mathbb{Q}} \, \mathbb{C}/\Lambda_f \, .
\end{align}
In other words, if $f \in \mathcal{M}_k(\Gamma)$ has algebraic residues, the associated class $C_f \in H^1(\Gamma,V_{k-2}(\mathbb{C}/\Lambda_f))$ can be identified with a linear combination of basis elements of $H^1(\Gamma,V_{k-2}(\mathbb{Q}))$ with coefficients in $\mathbb{C}/\Lambda_f$. We call these coefficients the periods of $f$. The conjecture above then states that magneticity is equivalent to the vanishing of periods. 
We now give examples which, among others, motivate the conjecture (for more examples see Appendix~\ref{sec:NumericalPeriods}).

% If $f$ has poles, then the value of $C_f(\gamma)$ will in general depend on the path from $\gamma^{-1}i\infty$ to $i\infty$, and the difference lies in the lattice $\frac{1}{(2\pi i)^{k/2-1}}\Lambda_f$ spanned by the residues of $f$:
% \beq
% \Lambda_f = \bigoplus_{\tau_p} (2\pi i)^{k/2}\,\Res_{\tau_p}f(\tau)\,\Z\,,
% \eeq
% where the sum runs over all poles $\tau_p$ of $f$. The normalisation is chosen such that, if $f$ is magnetic, Conjecture~\ref{conj:residues} implies that $\Lambda_f$ defines a lattice in $\overline{\Q}$. It is then natural to look at the class defined by $C_f$ in $H^1(\Gamma,\widetilde{V}_{k-2}^f)$,
% \beq
% \widetilde{V}_{k-2}^f := V_{k-2}/\big((2\pi i)^{1-k/2}\Lambda_f\big)
% \eeq
% is the vector space of homogeneous polynomials of degree $k-2$ with coefficients in $\C/(2\pi i)^{1-k/2}\Lambda_f$. 

% We have evaluated the $C_f(\gamma)$ numerically for all all known examples of strongly magnetic modular forms for various choices of $\gamma$. We find that in all cases the cohomology class $[C_f]$ can be represented in a very simple way, which we summarise in the following:

% Note that $Q(X,Y)$ is actually fixed for every choice of cusp $c$. Indeed, since $\delta(Q(X,Y))_{\gamma} = Q(X,Y)|_{\gamma}-Q(X,Y)$ and $C_f(\gamma)=0$ if $\gamma$ lies in the stabiliser $\Gamma_{c}$ of $c$ in $\Gamma$, we must have $Q(X,Y)|_{\gamma}=Q(X,Y)$ for every $\gamma\in\Gamma_c$.

\begin{example}
Let us compute the periods of $\phi$. Its residues were computed in Example~\ref{ex:phi_residues}.
Numerically we find that the periods of $\phi$ vanish. More precisely, for any $\gamma \in \Gamma_0(8)$ and any path of integration, the coefficients of the period polynomial $C_{f,i\infty}(\gamma)$ are in $(2\pi i)^2 \mathbb{Q}$. For example, integrating along the vertical line with $\text{Re } \tau = 1/8$ we find that 
\begin{align}
    C_{f,i\infty}(\gamma) \, = \, \frac{(2\pi i)^3}{2} \int_{1/8}^{i\infty} (X-\tau \, Y)^2 \, \phi(\tau) \, \rd \tau \, = \, (2\pi i)^2 \, (-\tfrac{1}{2} \, X^2 + \tfrac{1}{8} \, X \, Y)
\end{align}
with $\gamma=\mati{1}{0}{-8}{1}$.
\end{example}

\begin{example}
\label{exE}
Here we compute the periods of $E_4/j$ and $E_6/j$. Their residues were computed in Example~\ref{ex:SL2Ej}.
Numerically we find that the periods vanish. For example, integrating along the imaginary axis one finds that the period polynomials for $S = \mati{0}{-1}{1}{0}$ are
\beq\bsp
%\widehat{C}_{E_4(\tau)/j(\tau),i\infty}(S) &\,=
\frac{(2\pi i)^3}{2}\int_{0}^{i\infty} (X-\tau \, Y)^2 \, \frac{E_4(\tau)}{j(\tau)} \, \rd \tau 
&\,= (2\pi i)^2 \, \sqrt{-3} \, \frac{X \, Y}{576} +\omega_{E_4/j}\,(d^{(1)} Y^2)(S)\,,\\
%
%\widehat{C}_{E_6(\tau)/j(\tau),i\infty}(S) &\,=
\frac{(2\pi i)^5}{24}\int_{0}^{i\infty} (X-\tau \, Y)^4 \, \frac{E_6(\tau)}{j(\tau)} \, \rd \tau 
&\,= (2\pi i)^3 \, \frac{X^3 \, Y + X \, Y^3}{1536} +\omega_{E_6/j}\,(d^{(1)} Y^4)(S)\,,
\esp\eeq
with
\beq\bsp
\omega_{E_4/j} &\,:= \, \frac{(2\pi i)^3}{2}\int_0^{i\infty}\frac{E_4(\tau)}{j(\tau)} \, \rd\tau \, = \,  0.0610392510075\cdots\,,\\
\omega_{E_6/j} &\,:= \, \frac{(2\pi i)^5}{24}\int_0^{i\infty}\frac{E_6(\tau)}{j(\tau)} \, \rd\tau \, = \,  0.0876499825220\cdots\, .
\esp\eeq
In fact, in this case the periods vanish for trivial reasons. This is because these forms vanish at $i\infty$ and for $k=4,6$ every representative of a non-zero element of $H^1(\SL_2(\Z),V_{k-2}(\mathbb{Q}))$ cannot vanish on $\mati{1}{1}{0}{1}$ (equivalently, there are no cusp forms of level $1$ and weights $k=4,6$).

\end{example}

\begin{example}
Now consider the meromorphic modular form
\begin{align}
    \tilde{\phi} \, = \, \frac{(t^2+1)^4}{(t^2+2 \, t-1)^2 \, (t^2-2 \, t-1)^2} \, f \, \in \, \mathcal{M}_4(\Gamma_0(8)) \, ,
\end{align}
where $f$ and $t$ are given in eq.~\eqref{eq:Gamma0(8)}.
This is one of many similar examples in ref.~\cite{ThesisBoenisch}. The lattice generated by the residues is $\Lambda_{\tilde{\phi}} = (2\pi i)^2 \, \sqrt{-2} \, \mathbb{Q}$, but $\tilde{\phi}$ is not magnetic. The periods also do not vanish, but for any odd $n$, the periods of $T_n(\tilde{\phi}) - a_n \, \tilde{\phi}$ necessarily vanish, where $a_n$ is the $n^{\textrm{th}}$ Hecke eigenvalue (or equivalently the $n^{\textrm{th}}$ Fourier coefficient) of $f$. As predicted by our conjecture, the first terms in the Fourier expansions indeed suggest that the forms $T_n(\tilde{\phi}) - a_n \, \tilde{\phi}$ are magnetic.\footnote{This would be easy to understand if $\tilde{\phi}$ were the sum of an element in $\mathcal{S}_k(\Gamma_0(8))$ and a magnetic modular form, but this does not seem to be the case.} In other words, $\tilde{\phi}$ is a Hecke eigenform modulo magnetic modular forms. For example, we have
\begin{align}
    (T_3(\tilde{\phi})+4\, \tilde{\phi})(\tau) = \, &256 \, ( q + 2011260 \, q^3 + 2103414671870 \, q^5 + \cdots ) \, .
\end{align}

\end{example}
% !TEX root = main.tex

\section{Real-analytic modular forms with poles}
\label{sec:real_analytic}

In a recent series of papers~\cite{Brown:2018ut,brown_2020,Brown2017ACO}, Francis Brown has introduced a class of real-analytic modular forms.
%\footnote{Brown has only considered the case of real-analytic modular forms for the full modular group $\SL_2(\Z)$. There is no problem extending the definition to congruence subgroups, cf. also ref.~\cite{Drewitt}.} 
More precisely, a real-analytic modular form of weights $(r,s)$ for a finite-index subgroup $\Gamma$ of $\SL_2(\Z)$ is a real-analytic function $f:\HH\to \mathbb{C}$ such that for all $\mati{a}{b}{c}{d}\in\Gamma$ we have
\beq
f\left(\gamma\cdot\tau\right)  = f\left(\frac{a\tau+b}{c\tau+d}\right) = (c\tau+d)^r\,(c\bar{\tau}+d)^s\,f(\tau)\,.
\eeq
We now show that one can use certain meromorphic modular forms to construct real-analytic modular forms with poles. Our conjectures further suggest that magnetic modular forms are often suitable for this construction. 

Consider a meromorphic modular form $f \in \mathcal{M}_k(\Gamma)$ of weight $k \geq 2$ and again denote by $\Lambda_f \subset \mathbb{C}$ the $\mathbb{Q}$-vector space spanned by the residues $(2\pi i)^k \, \text{Res}_{\tau_0} \, \tau^m \, f(\tau)$ for $0\le m\le k-2$ and $\tau_0 \in \HH$. In the following we assume that $\Lambda_f$ is purely imaginary and that the cohomology class $C_f$ vanishes. For any choice of $\tau_0 \in \overline{\HH}$ at which $f$ is holomorphic we can consider the integral 
\beq\label{eq:F_def}
F_{\tau_0}(\tau,X,Y) \, := \, \frac{(2\pi i)^{k-1}}{(k-2)!} \, \int_{\tau_0}^{\tau}(X-\tau' Y)^{k-2} \, f(\tau') \, \rd\tau'\,. 
\eeq
This is only well-defined modulo $\Lambda_f$ and for every $\gamma \in \Gamma$ we have 
\begin{align*}
    F_{\tau_0}(\gamma \cdot \tau,X,Y)|\gamma \, = \, F_{\tau_0}(\tau,X,Y) + C_{f,\tau_0}(\gamma) \, .
\end{align*}
Now, since $C_f = 0$, there is a polynomial $P_{\tau_0} \in V_{k-2}(\mathbb{C})$ such that $C_{f,\tau_0} = d^{(1)}P_{\tau_0}$ and we find that $\tau \mapsto \mathcal{F}(\tau,X,Y) = \text{Re} (F_{\tau_0}(\tau,X,Y) - P_{\tau_0})$ defines a real-analytic function with poles $\HH \rightarrow V_{k-2}(\mathbb{R})$ that satisfies
\begin{align*}
    \cF(\gamma\cdot \tau,X,Y)|\gamma \, = \, \cF(\tau,X,Y) \, .
\end{align*}
To obtain real analytic modular forms with poles we can now apply Proposition~7.1 of ref.~\cite{Brown:2018ut} (see also Proposition 2.1 of ref.~\cite{Zemel2013OnQF}):
\begin{prop}\label{prop:prop}
Consider $\cF(\tau,X,Y) = \sum_{r+s=k-2}f_{r,s}(\tau)(X-\tau Y)^r(X-\bar{\tau} Y)^s$, where $f_{r,s}:\HH\to \C$ are real-analytic functions. Then $\cF(\gamma\cdot \tau,X,Y)|\gamma = \cF(\tau,X,Y)$ if and only if $f_{r,s}$ are real-analytic modular forms of weights $(r,s)$ for all $r+s=k-2$.
\end{prop}

Conjecture \ref{conj:periods} suggests that magnetic modular forms are suitable for the construction above, provided that the lattice $\Lambda_f$ is purely imaginary. 

\begin{example}
Consider the magnetic modular form $E_4/j$. The residues are purely imaginary (cf.~Example~\ref{ex:SL2Ej}) and the associated cocycle is a coboundary (cf.~Example~\ref{exE}). We obtain
\begin{align}
    \mathcal{F}(\tau,X,Y) \, &= \, \text{Re} \left(\frac{(2\pi i)^3}{2}\int_{i\infty}^\tau (X-\tau' \, Y)^2 \, \frac{E_4(\tau')}{j(\tau')} \, \rd \tau' -\omega_{E_4/j} \, Y^2 \right)
\end{align}
with $\omega_{E_4/j}$ from Example~\ref{exE} and it follows that the functions $f_{r,s}$ defined by 
\begin{align}
    \left( 
    \begin{array}{c}
         f_{2,0} \\
         f_{1,1} \\
         f_{0,2}
    \end{array}
    \right) (\tau) \, = \, \frac{1}{(\tau-\overline{\tau})^2} \left( 
    \begin{array}{ccc}
        \overline{\tau}^2 & -2 \, \overline{\tau} & 1  \\
         -2 \, \tau \, \overline{\tau} & 2 \, (\tau +\overline{\tau}) & -2 \\
        \tau^2 & -2 \, \tau & 1
    \end{array}
    \right) \, \left(
    \begin{array}{c}
         I_0  \\
         I_1  \\
         I_2-\omega_{E_4/j}
    \end{array}
    \right) (\tau)
\end{align}
in terms of 
\begin{align}
    I_n(\tau) \, = \, \text{Re} \left( \frac{(2\pi i)^3}{2}\int_{i\infty}^\tau \tau'^n \, \frac{E_4(\tau')}{j(\tau')} \, \rd \tau' \right)
\end{align}
are real-analytic modular forms with poles and weight $(r,s)$.
\end{example}

\begin{example}
Consider the magnetic modular form $E_6/j$. The residues are purely imaginary (cf.~Example~\ref{ex:SL2Ej}) and the associated cocycle is a coboundary (cf.~Example~\ref{exE}). We obtain
\begin{align}
    \mathcal{F}(\tau,X,Y) \, &= \, \text{Re} \left(\frac{(2\pi i)^5}{24}\int_{i\infty}^\tau (X-\tau' \, Y)^4 \, \frac{E_6(\tau')}{j(\tau')} \, \rd \tau' -\omega_{E_6/j} \, Y^4 \right)
\end{align}
with $\omega_{E_6/j}$ from Example~\ref{exE} and it follows that the functions $(f_{4,0},f_{3,1},...,f_{0,4})$ defined by 
{\tiny
    \begin{align}
    \frac{1}{(\tau-\overline{\tau})^4} \left( 
    \begin{array}{ccccc}
        \overline{\tau}^4 & -4 \, \overline{\tau}^3 & 6 \, \overline{\tau}^2 & -4 \, \overline{\tau} & 1  \\
        -4 \, \tau \, \overline{\tau}^3 & 12 \, \tau \, \overline{\tau}^2 + 4 \, \overline{\tau}^3 & -12 \, \tau \, \overline{\tau} - 12 \, \overline{\tau}^2 & 4 \, \tau +12 \, \overline{\tau} & -4  \\
        6 \, \tau^2 \, \overline{\tau}^2 & -12 \, \tau^2 \, \overline{\tau} -12 \, \tau \, \overline{\tau}^2 & 6 \, \tau^2 +24 \, \tau \, \overline{\tau} +6 \, \overline{\tau}^2 & -12 \, \tau -12 \, \overline{\tau} & 6  \\
        -4 \, \tau^3 \, \overline{\tau} & 4 \, \tau^3 + 12 \, \tau^2 \, \overline{\tau} & - 12 \, \tau^2 -12 \, \tau \, \overline{\tau} & 12 \, \tau + 4 \, \overline{\tau} & -4  \\
        \tau^4 & -4 \, \tau^3 & 6 \, \tau^2 & -4 \, \tau & 1
        \end{array}
    \right) \, \left(
    \begin{array}{c}
         I_0  \\
         I_1  \\
         I_2 \\
         I_3 \\
         I_4-\omega_{E_6/j}
    \end{array}
    \right) (\tau)
\end{align}
}%
in terms of 
\begin{align}
    I_n(\tau) \, = \, \text{Re} \left( \frac{(2\pi i)^5}{24}\int_{i\infty}^\tau \tau'^n \, \frac{E_6(\tau')}{j(\tau')} \, \rd \tau' \right)
\end{align}
are real-analytic modular forms with poles and weight $(r,s)$.
\end{example}

\section*{acknowledgements}
The authors are grateful to Herbert Gangl for organising the block seminar on `\emph{Picard-Fuchs eq's, irrationality proofs/measures, Ap\'ery motives and related topics}', where some of the ideas presented here were first discussed. The authors thank Matt Kerr and Wadim Zudilin for discussions. K.B. is supported by the International Max Planck Research School on Moduli Spaces of the Max Planck Institute for Mathematics in Bonn. This work was co-funded by the European Union through the ERC Consolidator Grant LoCoMotive 101043686. Views and opinions expressed are however those of the author(s) only and do not necessarily reflect those of the European Union or the European Research Council. Neither the European Union nor the granting authority can be held responsible for them.

\appendix
% !TEX root = main.tex

\section{Finding candidates for magnetic modular forms}
\label{app:fit}

In this appendix we explain a method for determining candidates of magnetic modular forms. For any finite subset $R \subset \overline{\HH}/\Gamma$ we denote by $\mathcal{M}_k^R(\Gamma)$ the subset of $\mathcal{M}_k(\Gamma)$ consisting of meromorphic modular forms which are holomorphic away from $R$. The quotient $\mathcal{M}_k^R(\Gamma) / D^{k-1}\mathcal{M}_{2-k}^R(\Gamma)$ is finite-dimensional and in practice one can often write down a basis. Magnetic modular forms in $\mathcal{M}_k^R(\Gamma)$ then just correspond to solutions of an (infinite-dimensional) linear system of congruence equations. To find candidates for magnetic modular forms, one can take finitely many equations from this system and compute solutions using modern computer algebra systems.

If $\Gamma$ has genus 0, writing down a basis of $\mathcal{M}_k^R(\Gamma) / D^{k-1}\mathcal{M}_{2-k}^R(\Gamma)$ becomes particularly easy. For simplicity we further restrict to $k > 2$ and $R$ consisting of all cusps and the orbit of some $\tau_0 \in \HH$. We can then use a result of refs.~\cite{matthes2021iterated,Broedel:2021zij}, which we state here only in the case where $\Gamma$ has no irregular cusps and no elliptic points, namely that representatives of generators of $\mathcal{M}_k^R(\Gamma) / D^{k-1}\mathcal{M}_{2-k}^R(\Gamma)$ are given by the following meromorphic modular forms of weight $k$ for $\Gamma$:
\begin{itemize}
\item forms with poles only in the orbit of $i\infty$ and of order at most $\dim S_k(\Gamma)$
\item forms with poles only in the orbit of $\tau_0$ and of order at most $k-1$
\end{itemize}

We applied this method to determine the candidates for magnetic modular forms of weights 4 and 6 for $\Gamma(2)$ and $\Gamma_1(6)$ presented in Appendices~\ref{sec:Gamma(2)} and~\ref{sec:Gamma1(6)}. Note that, while this procedure was able to determine those candidates of even weights, the congruence conditions did not seem to have solutions for odd weights, in agreement with our Conjecture~\ref{conj:pole}.

% !TEX root = main.tex

\section{New candidates for magnetic modular forms}
\label{app:candidates}

Here we present several new candidates for magnetic modular forms of small weights and depths for various congruence subgroups. We do not have a proof that these functions are magnetic, but we have checked that (at least) the first 500 Fourier coefficients have the correct divisibility properties. We therefore conjecture that these modular forms are indeed magnetic.

\subsection{Candidates for magnetic modular forms for $\boldsymbol{\Gamma_0(2)}$}
\label{sec:Gamma0(2)}
We start by presenting three candidates of magnetic modular forms for $\Gamma_0(2)$. All three candidates were found by extending the construction of ref.~\cite{Pogel:2022yat} to a hypergeometric family of K3 surfaces whose Picard-Fuchs operator is the symmetric square of an operator of degree two with monodromy group $\Gamma_0(2)$~\cite{epsilon_form}.
These three functions read
\beq\bsp\label{eq:Gamna0(2)_candidates}
C_4&\,:= \frac{A}{3(1+64 A)}\,\left(4\,\phi_0^2-E_4\right)\,,\\
%&\,\phantom{:} = q-56\, q^2+2076\, q^3-65984\, q^4+1941630\, q^5+\ldots\,,\\
C_{6a}&\,:= \frac{A\,(5+1408 A+20480 A^2)}{9(1-64 A)^3}\,\left(8\,\phi_0^3+E_6\right)\,, \\
%&\,\phantom{:}= 5\,q+2528\, q^2+547524\, q^3+87849984\, q^4+12091540750\, q^5+\ldots\,,\\
C_{6b} &\,:= \frac{A\,(1-64 A)}{9(1+64 A)^2}\,\left(8\,\phi_0^3+E_6\right)\,,
%&\, \phantom{:}= q-160\, q^2+9972\, q^3-447488\, q^4+17028150\, q^5+\ldots\,,\\
\esp\eeq
where $\phi_0$ is the unique normalised Eisenstein series of weight 2 for $\Gamma_0(2)$,
\beq
\phi_0(\tau)  :=  1+24\sum_{m=1}^\infty\frac{mq^m}{1+q^m} = 1+24q+24q^2+96q^3+\cdots\,,
\eeq
and $A$ is the following Hauptmodul for $\Gamma_0(2)$:
\beq
A(\tau) := \frac{\eta(2\tau)^{24}}{\eta(\tau)^{24}} = q +24q^2 + 300q^3+\cdots\,.
\eeq
%From the definitions it is clear that $C_4$, $C_{6a}$ and $C_{6b}$ are modular forms of weights 4 and 6 respectively. 
From the divisibility of the Fourier coefficients we see that $C_4$ is a candidate for a magnetic modular form of weight 4 and depth 1, and $C_{6a}$ and $C_{6b}$ are candidates of weight 6 and depth 2. 
%In all cases we have chosen a normalisation so that the number $b$ from Definition~\ref{def:magnetic} is equal to unity.

%%%%%%%%%%%%%%%%%%%%%%%%%%%%%%%%%%%%%%%%%%%%%%%%%%

\subsection{Candidates for magnetic modular forms for $\boldsymbol{\Gamma(2)}$}
\label{sec:Gamma(2)}
The following are candidates for magnetic modular forms of weight 4 and depth 1:
\beq\bsp
%D_{4a}(\tau) &\,:= \frac{(1-\lambda (\tau )) \lambda (\tau )^2 }{4\,(\lambda (\tau )-2)^2}\,\rho (\tau )^2\,,
%&\,\phantom{:}=q-56\, q^2+2 076\, q^3-65984 q^4+1941630 q^5+\ldots\,,\\
D_{4a} &\,:=\frac{4\,(1-\lambda)^2 \lambda }{(\lambda +1)^2}\,\rho^2\,,\\
%&\,\phantom{:}=q-56 q^2+2076 q^3-65984 q^4+1941630 q^5+\ldots\,,\\
D_{4b} &\,:= \frac{4\,(1-\lambda ) \lambda  }{(1-2 \lambda )^2}\,\rho^2\,,
%&\,\phantom{:}=\,,\\
\esp\eeq
where $\rho = \theta_3^4$ is an Eisenstein series of weight two for $\Gamma(2)$, $\lambda = \theta_2^4/\theta_3^4$ is a Hauptmodul for $\Gamma(2)$ and $\theta_3$ and $\theta_4$ are the usual theta functions.
We also have candidates for magnetic modular forms of weight 6 and depth 2:
\beq\bsp
%D_{4a}(\tau) &\,:= \frac{(1-\lambda (\tau )) \lambda (\tau )^2 }{4\,(\lambda (\tau )-2)^2}\,\rho (\tau )^2\,,
%&\,\phantom{:}=q-56\, q^2+2 076\, q^3-65984 q^4+1941630 q^5+\ldots\,,\\
D_{6a} &\,:=\frac{4 (1-\lambda)^2 \lambda (1+\lambda )\left(5 \lambda^4+68 \lambda^3-66 \lambda^2+68 \lambda+5\right)  }{\left(\lambda^2-6 \lambda+1\right)^3}\,\rho^3 \,,\\
%&\,\phantom{:}=q-56 q^2+2076 q^3-65984 q^4+1941630 q^5+\ldots\,,\\
D_{6b}&\,:= -\frac{4 (1-\lambda) \lambda (2 \lambda-1) \left(80 \lambda^4-160 \lambda^3+168 \lambda^2-88 \lambda+5\right)}{\left(4 \lambda^2-4 \lambda-1\right)^3}\,\rho^3  \,,\\
D_{6c} &\,:=
\frac{4 (1-\lambda)^2 \lambda \left(\lambda^2-6 \lambda+1\right) }{(\lambda+1)^3}\,\rho^3\, , \\
D_{6d} &\,:= -\frac{4 (1-\lambda) \lambda \left(4 \lambda^2-4 \lambda-1\right) }{(2 \lambda-1)^3}\,\rho^3\,.
%&\,\phantom{:}=\,,\\
\esp\eeq
The candidate $D_{6d}$ was found by applying the same method as in the $\Gamma_0(2)$ case, but to a family of K3 surfaces whose periods can be written as modular forms of $\Gamma(2)$~\cite{epsilon_form}.
The other candidates were found using the strategy from Appendix~\ref{app:fit}.

\subsection{Candidates for magnetic modular forms for $\boldsymbol{\Gamma_1(6)}$}
\label{sec:Gamma1(6)}
The following are candidates of weight 4 and depth 1 and weight 6 and depth 2 respectively:
\beq\bsp
H_{4a} &\, :=\frac{2  (\xi -9) (\xi -1) \xi ^2}{(\xi-3)^2}\,\aleph_1^4\,,\\
H_{4b}&\, :=\frac{(\xi -9) (\xi -1) \xi  \left(\xi ^2-10 \xi +9\right)}{(\xi +3)^2}
\aleph_1^4\,,\\
H_{6a} &\, :=\frac{4 (\xi -9) (\xi -1) \xi ^2 (\xi +3) \left(\xi ^2-12 \xi +9\right)}{(\xi -3)^3}\,\aleph_1^6\,,\\
H_{6b} &\, :=\frac{ (\xi -9)^2 (\xi -3) (\xi -1)^2 \xi  \left(\xi^2-42 \xi +9\right)}{(\xi +3)^3}\,\aleph_1^6\,,
\esp\eeq
where $\aleph_1$ is a modular form of weight 1 and $\xi$ is a Hauptmodul for $\Gamma_1(6)$:
\beq
\aleph_1(\tau) = \frac{\eta (\tau ) \eta (6 \tau )^6}{\eta (2 \tau )^2 \eta (3 \tau )^3}
\qquad\textrm{~~~~~and~~~~~}\qquad
\xi(\tau) = \frac{\eta (2 \tau )^8 \eta (3 \tau )^4}{\eta (\tau )^4 \eta (6 \tau )^8}\,.
\eeq
All four examples were found using the strategy from Appendix~\ref{app:fit}. By comparing Fourier expansions, we find that
\beq\label{eq:Hab}
H_{4a}(\tau) = 2\,H_{4b}(2\tau)\qquad\textrm{~~~and~~~}\qquad  H_{6a}(\tau) = 4\,H_{6b}(2\tau)\,.
\eeq
Note that the weight 6 example of ref.~\cite{Pogel:2022yat} has triple poles at $\xi(\tau)=\pm3$ and can be written as a linear combination of $H_{6a}$ and $H_{6b}$.

\section{Numerical evidence for our conjectures}
\label{app:evidence}

\subsection{Action by $\boldsymbol{\SL_2(}\protect\fakebold{\mathbb{Z}}\boldsymbol{)}$}
\label{sec:NumericalSL2Action}
\subsubsection{Examples for $\Gamma(2)$}
We have checked how $\SL_2(\Z)$ acts on the candidates for magnetic modular forms given in the previous section. For the candidates of level 2, we find that they span three-dimensional invariant subspaces inside $\cMmag_{k,d}(\Gamma(2))$. More precisely, if we act with $S=\left(\begin{smallmatrix} 0&1\\-1&0\end{smallmatrix}\right)$ and $T=\left(\begin{smallmatrix} 1&1\\0&1\end{smallmatrix}\right)$ (which generate $\SL_2(\Z)$), we find:
\beq\bsp\label{eq:C4_SL2}
C_{4|S} = \frac{1}{4}\,D_{4a}\,, \qquad & D_{4a|S} = \phantom{-}4C_4\,,\qquad D_{4b|S} = \phantom{-}D_{4b}\,,\\
C_{4|T} =C_4\,,\!\phantom{\frac{1}{4}}\phantom{\frac{1}{4}}\qquad &D_{4a|T} = -D_{4b}\,,\qquad D_{4b|T} = -D_{4a}\,,
\esp\eeq
\beq\bsp\label{eq:C6a_SL2}
C_{6a|S} = \frac{1}{8}\,D_{6a}\,,\qquad &D_{6a|S} = 8C_{6a}\,,\qquad D_{6b|S} = D_{6b}\,,\\
C_{6a|T} =C_{6a}\,,\phantom{ \frac{1}{8}}\qquad &D_{6a|T} = \phantom{8}D_{6b}\,,\qquad D_{6b|T} = D_{6a}\,,
\esp\eeq
\beq\bsp\label{eq:C6b_SL2}
C_{6b|S} = \frac{1}{8}\,D_{6c}\,,\qquad &D_{6c|S} = 8C_{6b}\,,\qquad  D_{6d|S} = D_{6d}\,,\\
C_{6b|T} =C_{6b}\,,\phantom{ \frac{1}{8}}\qquad &D_{6c|T} = \phantom{8}D_{6d}\,,\qquad D_{6d|T} = D_{6c}\,.
\esp\eeq

\subsubsection{Examples for $\Gamma_1(6)$} $\Gamma_1(6)$ has index 24 in $\SL_2(\Z)$. We choose the following coset representatives:
\beq\bsp
&\gamma_1=\left(
\begin{smallmatrix}
 1 & 0 \\
 0 & 1 \\
\end{smallmatrix}
\right)\,,\qquad
\gamma_2=\left(
\begin{smallmatrix}
 0 & -1 \\
 1 & 0 \\
\end{smallmatrix}
\right)\,,\qquad
\gamma_3=\left(
\begin{smallmatrix}
 1 & 0 \\
 1 & 1 \\
\end{smallmatrix}
\right)\,,\qquad
\gamma_4=\left(
\begin{smallmatrix}
 0 & -1 \\
 1 & 2 \\
\end{smallmatrix}
\right)\,,\qquad
\gamma_5=\left(
\begin{smallmatrix}
 0 & -1 \\
 1 & 3 \\
\end{smallmatrix}
\right)\,,\qquad\\
&\gamma_6=\left(
\begin{smallmatrix}
 0 & -1 \\
 1 & 4 \\
\end{smallmatrix}
\right)\,,\qquad
\gamma_7=\left(
\begin{smallmatrix}
 0 & -1 \\
 1 & 5 \\
\end{smallmatrix}
\right)\,,\qquad
\gamma_8=\left(
\begin{smallmatrix}
 1 & 0 \\
 2 & 1 \\
\end{smallmatrix}
\right)\,,\qquad
\gamma_9=\left(
\begin{smallmatrix}
 1 & 1 \\
 2 & 3 \\
\end{smallmatrix}
\right)\,,\qquad
\gamma_{10}=\left(
\begin{smallmatrix}
 1 & 2 \\
 2 & 5 \\
\end{smallmatrix}
\right)\,,\qquad\\
&\gamma_{11}=\left(
\begin{smallmatrix}
 1 & 0 \\
 3 & 1 \\
\end{smallmatrix}
\right)\,,\qquad
\gamma_{12}=\left(
\begin{smallmatrix}
 -1 & -1 \\
 3 & 2 \\
\end{smallmatrix}
\right)\,,\qquad
\gamma_{13}=\left(
\begin{smallmatrix}
 -1 & -2 \\
 6 & 11 \\
\end{smallmatrix}
\right)\,,\qquad
\gamma_i = \gamma_{13}\gamma_{i-13}\,,\,\,\,i>13\,.
\esp\eeq
When acting with each $\gamma_i$, $1<i\le 12$ on the candidates for magnetic modular forms for $\Gamma_1(6)$ from Appendix~\ref{sec:Gamma1(6)}, we obtain in each case 11 new candidates for magnetic modular forms for $\Gamma(6)$. The remaining 12 coset representatives do not produce any new candidates. In particular, we have
\beq
H_{4a}|\gamma_{13} = H_{4a}\qquad \textrm{~~~and~~~}\qquad H_{6a}|\gamma_{13} = H_{6a}\,.
\eeq

\subsection{Action by Atkin-Lehner involutions}
\label{sec:NumericalALAction}
\subsubsection{Examples for $\Gamma_0(2)$} We find that the candidates of modular forms for $\Gamma_0(2)$ from section~\ref{sec:Gamma0(2)} are left invariant by $W_2$.

\subsubsection{Examples for $\Gamma_1(6)$} We act with $W_2$ and $W_3$ on the candidates from section~\ref{sec:Gamma1(6)}. We find that $W_3$ leaves all four candidates invariant, while $W_2$ acts according to
\beq\bsp
W_2(H_{4a}) = \frac{1}{2}H_{4b}\,,\qquad W_2(H_{4b}) = 2 \, H_{4a}\,,\\
W_2(H_{6a}) = \frac{1}{2}H_{6b}\,,\qquad W_2(H_{6b}) = 2 \, H_{6a}\,.
\esp\eeq

%\kiliancomment{Temporary collection of some results: 
%\begin{itemize}
%    \item[-] $C_4,C_{6a},C_{6b}$ are invariant under $W_2$
%    \item[-] $W_2(H_{4a}) = \frac{1}{2}H_{4b}$ and hence also $W_2(H_{4b}) = 2 \, H_{4a}$
%    \item[-] $W_2(H_{6a}) = \frac{1}{2}H_{6b}$ and hence also $W_2(H_{6b}) = 2 \, H_{6a}$
%\item[-] $H_{4a},H_{4b},H_{6a},H_{6b}$ are invariant under $W_3$
%\end{itemize}
%}

\subsection{Poles and residues}
\label{sec:NumericalResidues}

\subsubsection{Examples for $\SL_2(\Z)$}
We now discuss some of the magnetic modular forms for $\SL_2(\Z)$ from refs.~\cite{magnetic2,magnetic3,magnetic4}. We only show results for a selection of the functions considered in those references, but we have checked that the conclusions remain unchanged for all examples presented there.
The functions $E_4/j$ and $E_6/j$ were already discussed in Example \ref{ex:SL2Ej} and here we consider in addition
\beq\bsp\label{eq:Fab}
F_{4} &\,:= \frac{E_4}{j-1728}\,,\\
F_{8a} &\,:= \frac{j-2^{10}\cdot 3}{j^2} \, E_4^2\,,\\
F_{8b} &\,:= \frac{\big(13 j^3 - 443556 j^2 + 
    1610452125 j - 98280 \cdot 15^6)}{(j + 15^3)^4} \, E_4^2\,.
\esp\eeq
The residues are 
\beq\bsp
\Res_{\sqrt{-1}}\,\matix{1\\\tau\\\tau^2}F_{4}(\tau) &\,= \, \frac{\sqrt{-1}}{(2\pi i)^2}\,\matix{1/432\\0\\1/432}\,, \\
\Res_{-\frac{1}{2}+\frac{\sqrt{-3}}{2}}(\tau)\,\matix{1\\\tau\\\tau^2\\\tau^3\\\tau^4\\\tau^5\\\tau^6} F_{8a} &\,= \frac{\sqrt{-3}}{(2\pi i)^4}\,\matix{5/27\\-5/54\\2/27\\-7/108\\2/27\\-5/54\\5/27} \, , \\
\Res_{-\frac{1}{2}+\frac{\sqrt{-7}}{2}}(\tau)\,\matix{1\\\tau\\\tau^2\\\tau^3\\\tau^4\\\tau^5\\\tau^6} F_{8b} &\,= \frac{\sqrt{-7}}{(2\pi i)^4}\,\matix{40/2401\\-20/2401\\24/2401\\-26/2401\\48/2401\\-80/2401\\320/2401}\,.
\esp\eeq

\subsubsection{Examples for $\Gamma_0(2)$}
Let us consider the candidates for magnetic modular forms for $\Gamma_0(2)$ from Appendix~\ref{sec:Gamma0(2)}. The residues are 
\beq\bsp
\Res_{-\frac{1}{2}+\frac{\sqrt{-1}}{2}} \, \matix{1\\\tau\\\tau^2} C_4(\tau) &\,= \frac{\sqrt{-1}}{(2\pi i)^2}\,\matix{-1/8\\1/16\\-1/16}\,,\\
\Res_{\frac{\sqrt{-2}}{2}} \, \matix{1\\\tau\\\tau^2\\\tau^3\\\tau^4}C_{6a}(\tau) &\,= \frac{1}{(2\pi i)^3}\,\matix{3/4\\0\\1/8\\0\\3/16}\,,\\
\Res_{-\frac{1}{2}+\frac{\sqrt{-1}}{2}} \, \matix{1\\\tau\\\tau^2\\\tau^3\\\tau^4}C_{6b}(\tau) &\,= \frac{1}{(2\pi i)^3}\,\matix{-3/2\\3/4\\-1/2\\3/8\\-3/8}\,.
\esp\eeq
Since the candidates for magnetic modular forms for $\Gamma(2)$ from Appendix~\ref{sec:Gamma(2)} are related to $C_4$, $C_{6a}$ and $C_{6b}$ via eqs.~\eqref{eq:C4_SL2}, \eqref{eq:C6a_SL2} and~\eqref{eq:C6b_SL2}, we do not discuss those examples separately.

\subsubsection{Examples for $\Gamma_1(6)$}
    Let us consider the candidates for magnetic modular forms for $\Gamma_1(6)$ from Appendix~\ref{sec:Gamma1(6)}. We only show the results for $H_{4b}$ and $H_{6b}$, because the other are related to those by eq.~\eqref{eq:Hab}. The residues are 
\beq\bsp
    \Res_{\frac{1}{2}+\frac{\sqrt{-3}}{6}}\,\matix{1\\\tau\\\tau^2}H_{4b}(\tau) &\,= \frac{\sqrt{-3}}{(2\pi i)^2}\,\matix{-2/3\\-1/3\\-2/9}\,,\\
    \Res_{\frac{1}{2}+\frac{\sqrt{-3}}{6}}\,\matix{1\\\tau\\\tau^2\\\tau^3\\\tau^4}H_{6b}(\tau) &\,= \frac{1}{(2\pi i)^3}\,\matix{-36\\-18\\-10\\-6\\-4}\,.
    \esp\eeq

\subsection{Periods}
\label{sec:NumericalPeriods}

\subsubsection{Examples for $\SL_2(\Z)$}
We now discuss the examples for the full modular group of refs.~\cite{magnetic2,magnetic3,magnetic4}. The functions $E_4/j$ and $E_6/j$ were already discussed in Example \ref{exE} and here we only show the results for the functions in eq.~\eqref{eq:Fab}, but we have checked that the same conclusions hold for all strongly modular forms discussed in refs.~\cite{magnetic2,magnetic3,magnetic4}. Integrating along the imaginary axis (and avoiding poles by passing on the right) we find that
\beq\bsp
C_{F_{4},i\infty}(S) &\,= (2\pi i)^2 \, \sqrt{-1} \, \frac{X^2-4 \, X \, Y + Y^2}{1728} + \omega_{F_{4}}\,(d^{(1)} Y^2)(S)\,,\\
C_{F_{8a},i\infty}(S) &\,=(2\pi i)^4 \, \sqrt{-3} \, \frac{- X^5 \, Y -X^3 \, Y^3- X \, Y^5}{11664} + \omega_{F_{8a}}\,(d^{(1)} Y^6)(S)\,,\\
C_{F_{8b},i\infty}(S) &\,= (2\pi i)^4 \, \sqrt{-7} \, \frac{-6 \, X^5 \, Y -5 \, X^3 \, Y^3-6 \, X \, Y^5}{43218}+ \omega_{F_{8b}}\,(d^{(1)} Y^6)(S)
\esp\eeq
with
\beq\bsp
\omega_{F_{4}} &\,:= \text{Re} \left( \frac{(2\pi i)^3}{2}\int_0^{i\infty}F_4(\tau) \, \rd\tau\right) \, = \,  -0.0424058145452\cdots \,,\\
\omega_{F_{8a}} &\,:= \frac{(2\pi i)^7}{720}\int_0^{i\infty}F_{8a}(\tau) \, \rd\tau = 0.1175032101096\cdots\,,\\
\omega_{F_{8b}} &\,:= \frac{(2\pi i)^7}{720}\int_0^{i\infty}F_{8b}(\tau) \, \rd\tau = 0.3506647091986\cdots\,.
\esp\eeq

\subsubsection{Examples for $\Gamma_0(2)$}
$\Gamma_0(2)$ is generated by $T$ and $R= \mati{1}{0}{2}{1}$. Integrating along the axis with $\text{Re} \, \tau = -\frac{1}{2}$ (and avoiding poles by passing on the right) we find that
\beq\bsp
C_{C_{4},i\infty}(R) &\,= (2\pi i)^2 \, \sqrt{-1} \, \frac{-X^2}{32} + \omega_{C_{4}}\,(d^{(1)} Y^2)(R)\,,\\
C_{C_{6a},i\infty}(R) &\,=(2\pi i)^3 \, \frac{-4 \, X^3 \, Y-6 \, X^2 \, Y^2-4 \, X \, Y^3 - Y^4}{128} + \omega_{C_{6a}}\,(d^{(1)} Y^4)(R)\,,\\
C_{C_{6b},i\infty}(R) &\,= (2\pi i)^3 \, \frac{-2 \, X^4 - 2 \, X^3 \, Y- X^2 \, Y^2}{64} + \omega_{C_{6b}}\,(d^{(1)} Y^4)(R)
\esp\eeq
with
\beq\bsp
\omega_{C_{4}} &\,:= \frac{1}{4} \, \text{Re} \left( \frac{(2\pi i)^3}{2}\int_{-\frac{1}{2}}^{i\infty}C_4(\tau) \, \rd\tau\right) \, = \,  0.2289913985443\cdots \,,\\
\omega_{C_{6a}} &\,:= \frac{1}{16} \, \frac{(2\pi i)^5}{24}\int_{-\frac{1}{2}}^{i\infty}C_{6a}(\tau) \, \rd\tau = -0.7888498426984\cdots\,,\\
\omega_{C_{6b}} &\,:= \frac{1}{16} \, \text{Re} \left(\frac{(2\pi i)^5}{24}\int_{-\frac{1}{2}}^{i\infty}C_{6b}(\tau) \, \rd\tau \right) = 0.2629499475661\cdots\,.
\esp\eeq
Since the candidates for magnetic modular forms for $\Gamma(2)$ from Appendix~\ref{sec:Gamma(2)} are related to $C_4$, $C_{6a}$ and $C_{6b}$ via eqs.~\eqref{eq:C4_SL2}, \eqref{eq:C6a_SL2} and~\eqref{eq:C6b_SL2}, we do not discuss those examples separately.

\subsubsection{Examples for $\Gamma_1(6)$}
    We here discuss the examples for the congruence subgroup $\Gamma_1(6)$ in Appendix~\ref{sec:Gamma1(6)}. We only show the results for $H_{4b}$ and $H_{6b}$, because the other are related to those by eq.~\eqref{eq:Hab}. $\Gamma_1(6)$ is generated by $T$, $\gamma_1 = \left(\begin{smallmatrix} -5&1\\-6&1\end{smallmatrix}\right)$ and $\gamma_2 = \left(\begin{smallmatrix} 7&-3\\12&-5\end{smallmatrix}\right)$.  
    Integrating along the axes with $\text{Re} \, \tau = \frac{1}{6}$ and $\text{Re} \, \tau = \frac{5}{12}$ we obtain
\beq\bsp
C_{H_{4b},i\infty}(\gamma_1) \,=& (2\pi i)^2 \, \sqrt{-3} \, \frac{X^2}{9} + \omega_{H_{4b}}\,(d^{(1)} Y^2)(\gamma_1)\,,\\
C_{H_{4b},i\infty}(\gamma_2) \,=& (2\pi i)^2 \, \sqrt{-3} \, \frac{13 \, X^2-10 \, X \, Y +2 \, Y^2}{9} + \omega_{H_{4b}}\,(d^{(1)} Y^2)(\gamma_2)\,,\\
C_{H_{6b},i\infty}(\gamma_1) \,=&(2\pi i)^3 \, \frac{13 \, X^4 -6 \, X^3 \, Y+ \, X^2 \, Y^2}{12} + \omega_{H_{6b}}\,(d^{(1)} Y^4)(\gamma_1)\,,\\
C_{H_{6b},i\infty}(\gamma_2) \,=&(2\pi i)^3 \, \frac{401 \, X^4 -652 \, X^3 \, Y+398 \, X^2 \, Y^2-108 \, X \, Y^3 +11 \, Y^4}{6} \\
&+ \omega_{H_{6b}}\,(d^{(1)} Y^4)(\gamma_2)
\esp\eeq

with

\beq\bsp
\omega_{H_{4b}} &\,:= \frac{1}{36} \, \text{Re} \left( \frac{(2\pi i)^3}{2}\int_{\frac{1}{6}}^{i\infty}H_{4b}(\tau) \, \rd\tau\right) \, = \, 0.3906512064482 \cdots \,,\\
\omega_{H_{6b}} &\,:= \frac{1}{1296} \, \text{Re} \left(\frac{(2\pi i)^5}{24}\int_{\frac{1}{6}}^{i\infty}H_{6b}(\tau) \, \rd\tau \right) = 0.4006856343865 \cdots\,.
\esp\eeq

\bibliography{bib}

\providecommand{\href}[2]{#2}\begingroup\raggedright\begin{thebibliography}{10}

\bibitem{Adams:2017ejb}
L.~Adams and S.~Weinzierl, {\it {Feynman integrals and iterated integrals of
  modular forms}},  {\em Commun. Num. Theor. Phys.} {\bf 12} (2018) 193--251,
  [\href{http://xxx.lanl.gov/abs/1704.08895}{{\tt 1704.08895}}].

\bibitem{magnetic2}
Y.~Li and M.~Neururer, {\it A magnetic modular form},  {\em Int. J. Number
  Theory} {\bf 15} (2019), no.~5 907--924.

\bibitem{magnetic1}
D.~J. Broadhurst and W.~Zudilin, {\it {A magnetic double integral}},  {\em J.
  Austral. Math. Soc.} {\bf 107} (2019), no.~1 9--25.

\bibitem{magnetic3}
V.~Pasol and W.~Zudilin, {\it Magnetic (quasi-)modular forms},  {\em Nagoya
  Math. J.} {\bf 248} (2022) 849--864.

\bibitem{ThesisBoenisch}
{Kilian B{\"o}nisch}, {\em Modularity of special motives of rank four
  associated with Calabi-Yau threefolds}.
\newblock PhD thesis, Rheinische Friedrich-Wilhelms-Universit{\"a}t Bonn,
  Sept., 2023.

\bibitem{FiberingOutCalabiYauMotives}
 K. Bönisch, V. Golyshev, and A. Klemm, \textit{Fibering out {C}alabi-{Y}au
  motives}. Work in progress.

\bibitem{Pogel:2022yat}
S.~P\"ogel, X.~Wang, and S.~Weinzierl, {\it {The three-loop equal-mass banana
  integral in \ensuremath{\varepsilon}-factorised form with meromorphic modular
  forms}},  {\em JHEP} {\bf 09} (2022) 062,
  [\href{http://xxx.lanl.gov/abs/2207.12893}{{\tt 2207.12893}}].

\bibitem{epsilon_form}
C.~Duhr, S.~Maggio, C.~Nega, L.~Tancredi, and F.~Wagner. {To appear}.

\bibitem{Bol}
G.~Bol, {\it {Invarianten linearer Differentialgleichungen}},  {\em Abh. Math.
  Sem. Univ. Hamburg} {\bf 16} (1949) 1--28.

\bibitem{EichlerAbelscheIntegrale}
M.~Eichler, {\it Eine verallgemeinerung der abelschen integrale},  {\em
  Mathematische Zeitschrift} {\bf 67} (1957) 267--298.

\bibitem{Brown:2018ut}
F.~Brown, {\it A class of non-holomorphic modular forms i},  {\em Research in
  the Mathematical Sciences} {\bf 5} (2018), no.~1 7.

\bibitem{brown_2020}
F.~Brown, {\it A class of nonholomorphic modular forms ii: Equivariant iterated
  eisenstein integrals},  {\em Forum of Mathematics, Sigma} {\bf 8} (2020) e31.

\bibitem{Brown2017ACO}
F.~Brown, {\it A class of non-holomorphic modular forms iii: real analytic cusp
  forms for $\mathrm{SL}_2(\mathbb{Z})$},  {\em Research in the Mathematical
  Sciences} {\bf 5} (2017) 1--36.

\bibitem{Zemel2013OnQF}
S.~Zemel, {\it On quasi-modular forms, almost holomorphic modular forms, and
  the vector-valued modular forms of shimura},  {\em The Ramanujan Journal}
  {\bf 37} (2013) 165--180.

\bibitem{matthes2021iterated}
N.~Matthes, {\it Iterated primitives of meromorphic quasimodular forms for
  {$\operatorname{SL}_2(\mathbb Z)$}},  {\em Trans. Amer. Math. Soc.} {\bf 375}
  (2022) 1443--1460, [\href{http://xxx.lanl.gov/abs/2101.11491}{{\tt
  2101.11491}}].

\bibitem{Broedel:2021zij}
J.~Broedel, C.~Duhr, and N.~Matthes, {\it {Meromorphic modular forms and the
  three-loop equal-mass banana integral}},  {\em JHEP} {\bf 02} (2022) 184,
  [\href{http://xxx.lanl.gov/abs/2109.15251}{{\tt 2109.15251}}].

\bibitem{magnetic4}
S.~L\"obrich and M.~Schwagenscheidt, {\it {Arithmetic properties of Fourier
  coefficients of meromorphic modular forms}},  {\em Algebra and Number Theory}
  {\bf 15} (2021), no.~9 2381--2401.

\end{thebibliography}\endgroup
\bibliographystyle{JHEP}

 \end{document}